\documentclass[11pt,leqno]{article}

\usepackage{amsmath,amsfonts,amscd,amssymb,theorem}

\long\def\comment#1\endcomment{}

\comment
\endcomment

\makeatletter
\begingroup
\gdef\th@dotted{\normalfont\itshape
  \def\@begintheorem##1##2{%
        \item[\hskip\labelsep \theorem@headerfont ##1\ ##2.]}%
\def\@opargbegintheorem##1##2##3{%
   \item[\hskip\labelsep \theorem@headerfont ##1\ ##2\ (##3).]}}
\endgroup
\makeatother

\theoremstyle{dotted}

\newtheorem{theorem}{Theorem}[section]
\newtheorem{lemma}[theorem]{Lemma}

\newtheorem{prop}[theorem]{Proposition}
\newtheorem{corr}[theorem]{Corollary}

\makeatletter
\begingroup
\gdef\th@upshape{\normalfont
  \def\@begintheorem##1##2{%
        \item[\hskip\labelsep \theorem@headerfont ##1\ ##2.]}%
\def\@opargbegintheorem##1##2##3{%
   \item[\hskip\labelsep \theorem@headerfont ##1\ ##2\ (##3).]}}
\endgroup
\makeatother

\theoremstyle{upshape}

\newtheorem{defn}[theorem]{Definition}
\newtheorem{remark}[theorem]{Remark}
\newtheorem{exa}[theorem]{Example}

\makeatletter
\renewcommand{\subsection}{\@startsection{subsection}{2}{0pt}{-3ex
plus -1ex minus -0.2ex}{-2mm plus -0pt minus
-2pt}{\normalfont\bfseries}} 
\renewcommand{\subsubsection}{\@startsection{subsubsection}{3}{0pt}{-3ex
plus -1ex minus -0.2ex}{-2mm plus -0pt minus
-2pt}{\normalfont\bfseries}} 
\makeatother

\makeatletter
\@addtoreset{equation}{section}
\makeatother

\newcommand{\cntrct}                
{\hspace{2pt}\raisebox{1pt}{\text{$\lrcorner$}}\hspace{2pt}}

\newcommand{\proof}[1][Proof.]{\smallskip\noindent{\em #1}}
\def\endproof{\hfill\ensuremath{\square}\par\medskip}

\renewcommand{\labelenumi}{{\normalfont(\roman{enumi})}}

\def\eqref#1{\thetag{\ref{#1}}}

\let\latexref=\ref
\def\ref#1{{\normalfont{\latexref{#1}}}}

\newcommand{\wt}{\widetilde}

\newcommand{\dg}{\dagger}

\newcommand{\idot}{{\:\raisebox{1pt}{\text{\circle*{1.5}}}}}
\newcommand{\hdot}{{\:\raisebox{3pt}{\text{\circle*{1.5}}}}}
\newcommand{\C}{\mathcal{C}}
\newcommand{\ppt}{{\operatorname{\sf pt}}}

\newcommand{\E}{\mathcal{E}}

\newcommand{\Hhom}{\operatorname{\mathcal{H}{\it om}}}
\newcommand{\Hom}{\operatorname{Hom}}
\newcommand{\Fun}{\operatorname{Fun}}
\newcommand{\Ext}{\operatorname{Ext}}

\newcommand{\V}{\operatorname{\sf V}}

\newcommand{\D}{\mathcal{D}}

\newcommand{\eps}{\epsilon}

\newcommand{\A}{\mathcal{A}}

\newcommand{\colim}{\operatorname{\sf colim}}
\renewcommand{\lim}{\operatorname{\sf lim}}

\newcommand{\Sets}{\operatorname{Sets}}
\newcommand{\Cat}{\operatorname{Cat}}

\newcommand{\id}{\operatorname{\sf id}}

\newcommand{\Y}{\operatorname{\sf Y}}

\newcommand{\Pos}{\operatorname{Pos}}
\newcommand{\Unf}{\operatorname{K}}
\newcommand{\unf}{\operatorname{k}}

\def\phi{\varphi}

\newcommand{\copr}{\sqcup}

\newcommand{\Sec}{\operatorname{Sec}}
\newcommand{\Tot}{\natural}

\newcommand{\Iso}{{\star}}

\newcommand{\bC}{\overline{\C}}

\newcommand{\dm}{\diamond}

\newcommand{\wM}{\widetilde{M}}

\newcommand{\Ex}{\operatorname{\mathcal{E}\mathit{x}}}

\newcommand{\Spl}{\operatorname{\mathcal{S}\mathit{p}\mathit{l}}}

\newcommand{\tw}{\operatorname{\sf tw}}

\newcommand{\fFun}{\operatorname{{\mathcal F}{\mathit u}{\mathit n}}}

\newcommand{\Ho}{\operatorname{Ho}}
\newcommand{\Hho}{\operatorname{\mathcal{H}\mathit{o}}}

\newcommand{\red}{\operatorname{\sf red}}

\title{On $2$-categories of extensions}

\author{D. Kaledin\thanks{Partially supported by the Ministry of
    Science and Higher Education of the Russian Federation
    (agreement no. 075-15-2025-303), in Steklov International
    Mathematical Center, and by the HSE University Basic Research
    Program.}}

\begin{document}

\maketitle

\tableofcontents

\section*{Introduction.}

It is well-known by now that a category obtained by localization
should not be treated simply as a category --- this forgets some
essential information informally known as
``enhancement''. Historically, this deficiency was first observed
when working with derived categories of abelian categories. A
barebone version of enhancement for such categories is given by the
triangulated structure of \cite{verd} but this is still not enough
for many pracical purposes: gluing does not work, functor categories
are not triangulated, and so on (for a slightly more detailed list
of problems, see e.g. \cite[Introduction]{kazh}).

The purpose of this short note is to address another very simple and
basic problem with unenhanced triangulated categories, and to show
how it can be resolved in the context of the general enhancement
formalism recently developed in \cite{ka.big}, \cite{ka.small}.

Namely, whatever notion of enhancement one uses, at the very least,
one expects, for any two objects $M,N \in \D$ in an enhanced
category $\D$, to have a homotopy type of maps $\Hhom(M,N)$ such
that the usual set of maps $\Hom(M,N)$ is the set
$\pi_0(\Hhom(M,N))$ of its connected components. If $\D$ is
triangulated --- for example, if it is the derived category of an
abelian category --- then one can recover higher homotopy groups as
$\pi_i(\Hhom(M,N)) = \Hom(M,N[-i])$. In particular, if both $M$ and
$N$ lie in the heart of a $t$-structure on $\D$, then all these
higher homotopy groups vanish, so that $\Hhom(M,N)$ is discrete
(abelian categories do not need an enhancement). At the next level
of difficulty, still fixing a $t$-structure, one can consider the
subcategory $\D_{[0,1]} \subset \D$ of objects sitting in
homological degrees $0$ and $1$. In this case, for any $M,N \in
\D_{[0,1]}$, the homotopy type $\Hhom(M,N)$ is $1$-truncated ---
that is, $\pi_i(\Hhom(M,N))=0$ for $i \geq 2$ --- and $1$-truncated
homotopy types are very efficiently described by groupoids. Thus in
this case, we expect to have groupoids of morphisms in $\D_{[0,1]}$,
so that it should have a natural structure of a $2$-category. We
call it the {\em $2$-category of extensions}.

Somewhat surprisingly, already at this point the triangulated
formalism does not deliver: it is not possible to recover the
$2$-category of extensions by only using a triangulated structure on
$\D$.

What we do in the paper is two things, then. Firstly, we restrict
ourselves to the case when $\D=\D(\A)$ is the derived category of an
abelian category $\A$, and we construct the $2$-category of
extensions corresponding to $\D_{[0,1]}(\A)$ by hand, that is, using
explicit chain complexes of objects in $\A$. Secondly, we show how
to do a very general construction using the enhancement formalism of
\cite{ka.big}, and we prove that the two constructions agree.

\bigskip

The paper is organized as follows. Section~\ref{cat.sec} contains
the necessary preliminaries on category theory, Grothendieck
construction and suchlike, mostly to fix terminology and notation
(both are the same as in \cite{ka.big} and \cite{ka.small}). We also
describe a convenient packaging of various $2$-categorical notions
based on the Grothendieck constrution, as in
\cite{ka.adj}. Section~\ref{ab.sec} is concerned with the classical
story: we recall some relevant facts on abelian categories and chain
complexes, especially complexes of length $2$, we construct our
$2$-category of extensions, and we give a complete description of
its categories of morphisms. Finally, in Section~\ref{der.sec}, we
re-do the same story using derived categories enhanced in the sense
of \cite{ka.big}. For simplicity, we assume that $\A$ has enough
injectives --- e.g. it could be a Grothendieck abelian category. We
recall very briefly the basics of the enhancement formalism and the
relevant general facts that we use as a black box (the main one is
the existence of semicartesian products). We then construct the
$2$-category using this formalism, and we prove that this simple
two-line construction indeed recovers the classical one.

\subsection*{Acknowledgements.} This paper arose as a side note to a
long project on square-zero extensions of abelian categories, joint
with W. Lowen and M. Booth; it is a pleasure to thank both Wendy and
Matt for many interesting and productive discussions. It is also a
pleasure to thank my colleagues at Steklov Institute, especially
Sasha Efimov and Anton Fonarev.

\section{Generalities on categories and $2$-categories.}\label{cat.sec}

\subsection{Categories and functors.}

We use the same notation and terminology as in \cite[Chapter
  1]{ka.big}. In particular, for any category $\C$, we denote the
opposite category by $\C^o$, and for any functor $\gamma:\C_0 \to
\C_1$, the opposite functor is $\gamma^o:\C_0^o \to \C_1^o$. For any
object $c \in \C$, the {\em left comma-fiber} $\C' /_\gamma c$ of a
functor $\gamma:\C' \to \C$ is the category of pairs $\langle
c',\alpha \rangle$, $c' \in \C$, $\alpha:\gamma(c') \to c$ a
morphism, and dually, {\em the right comma-fiber} is $c
\setminus_\gamma \C' = ({\C'}^o /_{\gamma^o} c)^o$. We drop $\gamma$
from notation and write $\C' / c$, $c \setminus \C'$ when there is
no danger of confusion. We denote by $\sigma(c):\C' / c \to \C'$
resp.\ $\tau(c):c \setminus \C' \to \C'$ the forgetful functor
sending $\langle c',\alpha \rangle$ to $c'$. The {\em fiber} $\C_c$
is the full subcategory $\C_c \subset \C / c$ of pairs $\langle
c',\alpha \rangle$ with invertible $\alpha$. A {\em section} of a
functor $\gamma:\C \to I$ is a functor $\sigma:I \to \C$ equipped
with an isomorphism $\gamma \circ \sigma \cong \id$. We denote by
$\ppt$ the point category --- one object, one morphism --- and for
any object $c \in \C$ in a category $\C$, we let $\eps_c:\ppt \to
\C$ be the embedding onto $c$.

A category is {\em essentially small} if it is equivalent to a small
category, and a functor $\gamma:\C' \to \C$ is {\em small} if
$\gamma^{-1}(\C_0)$ is essentially small for any essentially small
full subcategory $\C_0 \subset \C$. A category $\C$ is a {\em
  groupoid} if all its morphisms are invertible, and for any
category $\C$, the {\em isomorphism groupoid} $\C_{\Iso} \subset \C$
has the same objects as $\C$, and those morphisms between them that
are invertible. A category is {\em discrete} if all its morphisms
are identity maps. For any groupoid $\C$, the objects of the
discrete category $\pi_0(\C)$ are isomorphism classes of objects in
$\C$; we have the tautological functor $\C \to \pi_0(\C)$ that
admits a section $\pi_0(\C) \to \C$, unique up to an isomorphism. A
groupoid $\C$ is {\em connected} if $\pi_0(\C) = \ppt$. More
generally, a category $\C$ is {\em connected} if for any
decomposition $\C = \C_0 \copr \C_1$, exactly one of the components
$\C_l$, $l=0,1$ is empty. For any categories $\C$, $\C_0$, $\C_1$
equipped with functors $\gamma_l:\C_l \to \C$, $l=0,1$, we let $\C_0
\times_{\C} \C_1$ be the category of triples $\langle c_0,c_1,\alpha
\rangle$, $c_l \in \C_l$, $l=0,1$, $\alpha:\gamma_0(c_0) \to
\gamma_1(c_1)$ an isomorphism, and we denote $\C_0 \times_\C \C_1 =
\gamma_0^*\C_1$ when we want to emphasize the dependence on
$\gamma_0$.

Given a partially ordered set $J$, we treat it as a small category
in the usual way --- objects are elements $j \in J$, morphisms $j
\to j'$ are order relations $j \leq j'$. A subset $J' \subset J$ is
equipped with induced order unless indicated otherwise. We denote by
$\V$ the partially ordered set $\V = \{o,0,1\}$ with order relations
$o \leq 0,1$. Functors $c:\V \to \C$ to some category $\C$
correspond to diagrams
$$
\begin{CD}
c(0) @<<< c(o) @>>> c(1)
\end{CD}
$$
in $\C$, and dually for functors $\V^o \to \C$.

For any category $\E$ and any essentially small category $I$, we
denote by $\Fun(I,\E)$ the category of functors $I \to \E$; we
shorten $\Fun(I^o,\E)$ to $I^o\E$. We denote the category of
sections of a functor $\C \to I$ by $\Sec(I,\C)$. For any functor
$\gamma:I' \to I$ between essentially small categories, we denote by
$\gamma^*:\Fun(I,\E) \to \Fun(I',\E)$ the pullback functor given by
precomposition with $\gamma$, and we let
$\gamma_!,\gamma_*:\Fun(I',\E) \to \Fun(I,\E)$ be its left
resp.\ right-adjoint {\em Kan extension functors}, whenever they
exists. Left resp.\ right Kan extension with respect to the
projection $I \to \ppt$ to the point category $\ppt$ are the colimit
$\colim_I$ resp. the limit $\lim_I$. Pullbacks are limits over
$\V^o$, and pushouts are colimits over $\V$. A functor $\gamma:I'
\to I$ is {\em cofinal} if $i \setminus I'$ is connected for any $i
\in I$, and {\em final} if $\gamma^o$ is cofinal. A functor $\gamma$
that admits a left resp.\ right-adjoint is cofinal resp.\ final. For
any cofinal functor $\gamma:I' \to I$ between essentially small
categories, and any functor $E:I \to \E$ to some target category
$\E$, $\colim_IE$ exists iff so does $\colim_{I'}\gamma^*E$, and the
natural map $\colim_{I'}\gamma^*E \to \colim_IE$ is an isomorphism,
and dually for limits and final functors.

A category $\E$ is {\em complete} resp.\ {\em cocomplete} if
$\lim_IE$ resp. $\lim_IE$ exists for any functor $E:I \to \E$ from
an essentially small $I$. The category $\Sets$ of sets is complete
and cocomplete, and we have a fully faithful Yoneda embedding $\Y:I
\to I^o\Sets$. If $\E$ is complete resp.\ cocomplete, then right
resp.\ left Kan extensions exist for all $\E$-valued functors, and
are explicitly given by
\begin{equation}\label{kan.eq}
\gamma_!E(i) \cong \colim_{I' / i}\sigma(i)^*E, \qquad \gamma_*E(i)
\cong \lim_{i \setminus I'}\tau(i)^*E,
\end{equation}
for any functor $\gamma:I' \to I$ between essentially small
categories, object $i \in I$, and functor $E:I' \to \E$. This helps
in computations since one can often replace $I' / i$ resp.\ $i
\setminus I'$ with a cofinal resp.\ final subcategory. For example,
if $\gamma$ is fully faithful, then $I' / \gamma(i')$
resp.\ $\gamma(i') \setminus I'$ has a terminal resp.\ initial
object for any $i' \in I'$, and then $\gamma^* \circ \gamma_!E \cong
\gamma^* \circ \gamma_*E \cong E$ by \eqref{kan.eq}. Even if $\E$ is
not complete resp.\ cocomplete, but the limits or colimits in the
right-hand side of \eqref{kan.eq} exist, then so does the
corresponding Kan extension, and \eqref{kan.eq} still holds.

\subsection{Fibrations and cofibrations.}

We also assume known the ``Grothen\-dieck construction'' --- that
is, the machinery of fibrations and cofibrations of categories and
cartesian and cocartesian functors of \cite{SGA1}, see
e.g.\ \cite[Section 1.3]{ka.big}. Informally, fibrations
resp.\ cofibrations $\C \to I$ correspond to ``pseudofunctors'' from
$I^o$ resp.\ $I$ to the category of categories, with fibers $\C_i$,
$i \in I$ and transition functors $f^*:\C_{i'} \to \C_i$
resp.\ $f_!:\C_i \to \C_{i'}$ for any morphism $f:i \to i'$ in
$I$. For any fibration $\gamma:\C \to I$, the opposite functor
$\gamma^o:\C^o \to I^o$ is a cofibration, and vice versa. For any
cofibration $\E^\hdot \to \E$ with small fibers, and any target
category $\C$, the {\em relative functor category}
$\Fun(\E^\hdot|\E,\C)$ is the fibration $\Fun(\E^\hdot|\E,\C) \to
\E$ with fibers $\Fun(\E^\hdot_e,\C)$, $e \in \E$, and transition
functors $(f_!)^*$. A fibration $\C \to I$ is a {\em family of
  groupoids} if all its fibers $\C_i$, $i \in I$ are groupoids, or
equivalently, all maps in $\C$ are cartesian over $I$. For any
fibration $\C \to I$, we have a family of groupoids $\C_\flat
\subset \C$ with the same objects as $\C$, and those maps between
them that are cartesian over $I$ (if $I$ is the point category
$\ppt$, then $\C_\flat = \C_\Iso$ is the isomorphism groupoid of the
category $\C$). A fibration is {\em discrete} if its fibers are
discrete categories. For any family of groupoids $\C \to I$, we have
a canonical factorization
\begin{equation}\label{pi.0.fact}
\begin{CD}
\C @>{a}>> \pi_0(\C|I) @>{b}>> I,
\end{CD}
\end{equation}
where $b$ is the discrete fibration with fibers $\pi_0(\C_i)$, $i
\in I$, and $a$ is a family of groupoids with connected
fibers.

Discrete fibrations with small fibers correspond to honest functors
from $I^o$ to $\Sets$; explicitly, a functor $X:I^o \to \Sets$
corresponds to the forgetful fibration $IX \to I$, where the {\em
  category of elements} $IX$ has objects $\langle i,x \rangle$, $i
\in I$, $x \in X(i)$, with morphisms $\langle i,x \rangle \to
\langle i',x' \rangle$ given by morphisms $f:i \to i'$ such that
$X(f)(x')=x$. A morphism $X \to X'$ defines a functor $IX \to
IX'$. If we have a functor $X_{01}:I^o \to \Sets$ and two
subfunctors $X_0,X_1 \subset X_{01}$ such that $X_{01} = X_0 \cup
X_1$, with the intersection $X = X_0 \cap X_1$, then $X_{01} \cong
X_0 \copr_{X_{01}} X_1$ as soon as $I$ is small, so that the
cofibered coproduct is well-defined, and the comparison functor
\begin{equation}\label{IX.eq}
(IX_{01})^o\E \to (IX_0)^o\E \times_{(IX)^o\E} (IX_1)^o\E
\end{equation}
is an equivalence for any target category $\E$.

\begin{exa}\label{tw.exa}
For any category $I$, let $Y:I^o \times I \to \Sets$ be the
$\Hom$-pairing. Then the {\em twisted arrow category} $\tw(I)$ is
the category of elements $\tw(I) = (I \times I^o)Y$, with its
discrete fibration $\tw(I) \to I \times I^o$. The functor
$\sigma:\tw(I) \to I \times I^o \to I$ is a fibration with fibers
$\tw(I)_i \cong (I / i)^o$, and the functor $\tau:\tw(I) \to I
\times I^o \to I^o$ is a fibration with fibers $\tw(I)_i \cong i
\setminus I$.
\end{exa}

For any functor $\gamma:I' \to I$ and fibration resp.\ cofibration
$\C \to I$, the induced functor $\gamma^*\C \to I'$ is a fibration
resp.\ cofibration. We also use a right Kan extension operation for
fibrations of e.g.\ \cite[Section 1.4]{ka.big}; explicitly, for any
functor $\gamma:I' \to I$ and fibration $\C' \to I'$, the fibration
$\gamma_*\C' \to I$ has fibers
\begin{equation}\label{2kan.eq}
  (\gamma_*\C')_i \cong \Sec^\Tot(I' / i,\C'),
\end{equation}
where $\Sec^\Tot$ stands for the category of sections $I' / i \to
\C'$ cartesian over $I'$. We have the canonical adjunction functor
$\gamma^*\gamma_*\C' \to \C'$ cartesian over $I'$, and dually, for
any fibration $\C \to I$, we have the canonical adjunction functor
$\C \to \gamma_*\gamma^*\C$ cartesian over $I$.

\subsection{Nerves and $2$-categories.}\label{2cat.subs}

For $2$-categories, we use the packaging described in \cite{ka.adj}
and based on the Grothendieck construction. For any integer $n \geq
0$, we let $[n]$ be the ordinal $\{0,\dots,n\}$ with the usual
order, and we let $\Delta$ be the category of these ordinals $[n]$,
$n \geq 0$ and order-preserving maps between them. A {\em simplicial
  set} is a functor $X:\Delta^o \to \Sets$, and its {\em category of
  simplices} is the category of elements $\Delta X$. For any $n \geq
m \geq 0$, we denote by $s,t:[m] \to [n]$ the embedding identifying
$[m]$ with an initial resp.\ terminal segment of $[m]$. For any $n
\geq l \geq 0$, we then have a commutative square
\begin{equation}\label{seg.sq}
\begin{CD}
[0] @>{t}>> [l]\\
@V{s}VV @VV{s}V\\
[n-l] @>{t}>> [n]
\end{CD}
\end{equation}
in $\Delta$ that happens to be cocartesian. A functor $E:\Delta^o
\to \E$ to some category $\E$ {\em satisfies the Segal condition} if
it sends commutative squares \eqref{seg.sq} to cartesian squares in
$\E$. A {\em Segal category} is a fibration $\C \to \Delta$ such
that for any square \eqref{seg.sq}, the corresponding functor
$\C_{[n]} \to \C_{[l]} \times_{\C_{[0]}} \C_{[n-l]}$ is an
equivalence. As a trivial example, if we let $\eps=\eps_{[0]}:\ppt
\to \Delta$ be the embedding onto $[0]$, then $\eps_*\E \to \Delta$
is a Segal category for any category $\E$ --- indeed, by
\eqref{2kan.eq}, we have
\begin{equation}\label{eps.eq}
  \eps_*\E \cong \Fun([0] \setminus \Delta|\Delta,\E),
\end{equation}
where $[0] \setminus \Delta \to \Delta$ is the tautological discrete
cofibration, so that for any $[n] \in \Delta$, we have
$(\eps_*\E)_{[n]} \cong \E^{n+1}$. A {\em $2$-category} is a Segal
category $\C$ such that $\C_{[0]}$ is discrete. For any Segal
category $\C$, with some fiber $\C_{[0]}$ and its isomorphism
groupoid $\C_{[0],\Iso}$, we can choose a splitting
$\pi_0(\C_{[0],\Iso}) \to \C_{[0],\Iso}$ of the canonical projection
$\C_{[0],\Iso} \to \pi_0(\C_{[0]},\Iso)$, and construct the {\em
  reduction} $\C^{\red}$ by the cartesian square
\begin{equation}\label{red.sq}
\begin{CD}
\C^{\red} @>>> \C\\
@VVV @VV{\alpha}V\\
\eps_*\pi_0(\C_{[0],\Iso}) @>>> \eps_*\C_{[0]},
\end{CD}
\end{equation}
where $\alpha$ is the adjunction functor. Then $\C^{\red}$ is a
$2$-category. Informally, objects in a $2$-category $\C$ are objects
$c \in \C_{[0]}$, and the transition functors $s^*,t^*:\C_{[1]} \to
\C_{[0]}$ define a decomposition
\begin{equation}\label{C.cc}
\C_{[1]} \cong \coprod_{c,c' \in \C_{[0]}} \C(c,c'),
\end{equation}
where $\C(c,c')$ is the category of morphisms from $c$ to $c'$, with
$s^*$ resp.\ $t^*$ sending a morphism to its source resp.\ target.

\begin{exa}\label{repl.exa}
Let $\Delta^\hdot$ be the category of pairs $\langle [n],l\rangle$,
$[n] \in \Delta$, $l \in [n]$, with maps $\langle [n],l \rangle \to
\langle [n'],l' \rangle$ given by maps $f':[n] \to [n']$ such that
$f(l) \leq l'$. Then the forgetful functor $\Delta^\hdot \to \Delta$
is a cofibration with fibers $\Delta^\hdot_{[n]} \cong [n]$, and for
any category $I$, $\Delta^\Tot I = \Fun(\Delta^\hdot|\Delta,I)$ is a
Segal category, with fibers $(\Delta^\Tot I)_{[n]} \cong
\Fun([n],I)$. Its reduction $\Delta I = (\Delta^\Tot I)^{\red}$ of
\eqref{red.sq} is the {\em simplicial replacement} of $I$. Its
objects are the same as objects in $I$, and categories of morphisms
$\Delta I(i,i')$ of \eqref{C.cc} are discrete categories
corresponding to the $\Hom$-sets in $I$. In the square
\eqref{red.sq}, we have $(\Delta^\Tot I)_{[0]} \cong I$, and in
terms of \eqref{eps.eq}, the adjunction functor is induced by the
embedding $[0] \setminus \Delta \to \Delta^\hdot$.
\end{exa}

The simplicial replacement of a small category is the category of
simplices of a simplicial set. Namely, define the nerve functor
$N:\Cat \to \Delta^o\Sets$ from the category $\Cat$ of small
categories by
\begin{equation}\label{N.eq}
N = \phi_!\Y,
\end{equation}
where $\phi:\Delta \to \Cat$ is the standard embedding sending $[n]$
to itself treated as a small category, and $\Y:\Delta \to
\Delta^o\Sets$ is the Yoneda embedding. Then by \eqref{kan.eq}, we
have $\Delta I = \Delta N(I)$ for any $I \in \Cat$. Explicitly,
objects in $\Delta I$ are pairs $\langle [n],i \rangle$, $[n] \in
\Delta$, $i:[n] \to I$ a functor. If we take $I = [m]$ for some $[m]
\in \Delta$, then $\Delta [m] \cong \Delta / [m]$. The nerve functor
\eqref{N.eq} is fully faithful, and $X \in \Delta^o\Sets$ is in its
essential image iff it satisfies the Segal condition (or
equivalently, iff $\Delta X \to \Delta$ is a Segal category). Note,
however, that $N$ does not preserve coproducts. In particular, for
any cocartesian square \eqref{seg.sq} in $\Delta$, the maps
$N(s):N([l]) \to N([n])$, $N(t):N([n-l]) \to N([n])$ are injective,
we have $N([0]) \cong s(N([l])) \cap N([n-l])$, and \eqref{seg.sq}
induces an injective map
\begin{equation}\label{seg.eq}
b^l_n:N([l]) \copr_{N([0])} N([n-l]) \to N([n])
\end{equation}
whose image consists of arrows $[m] \to [n]$ that factors through $s:[l]
\to [n]$ or $t:[n-l] \to [m]$ (or both). Thus \eqref{seg.eq} is only
surjective if $l=0$ or $l=n$.

For a general $2$-category $\C$, the corresponding family of
groupoids $\C_\flat \to \Delta$ is also a $2$-category, with the
same objects as $\C$ and morphism categories $\C_\flat(c,c') =
\C(c,c')_\Iso$ given by isomorphism groupoids of the morphism
categories $\C(c,c')$. We can then consider the decomposition
\eqref{pi.0.fact}, and the discrete fibration
$\pi_0(\C_\flat|\Delta) \to \Delta$ is also a $2$-category with the
same objects, and morphism categories $\pi_0(\C(c,c')_\Iso)$. Say
that a $2$-category $\C$ is {\em bounded} if for any $c,c' \in
\C_{[0]}$, the category $\C(c,c')$ is essentially small. Then for a
bounded $2$-category $\C$, $\pi_0(\C_\flat|\Delta) \cong \Delta \bC$
for a unique category $\bC$ called the {\em truncation} of the
$2$-category $\C$. Its objects are objects of $\C$, and morphisms
are isomorphism classes of $1$-morphisms in $\C$.

A {\em Segal functor} between Segal categories $\C$, $\C'$ is a
functor $\C \to \C'$ cartesian over $\Delta$, and a $2$-functor
between $2$-categories $\C$, $\C'$ is a Segal functor $\C \to
\C'$. A Segal functor $\gamma:\C \to \C'$ is {\em $2$-fully
  faithful} if the square
\begin{equation}\label{ff.sq}
\begin{CD}
\C @>{\gamma}>> \C'\\
@V{\alpha}VV @VV{\alpha}V\\
\eps_*\C_{[0]} @>{\eps_*(\gamma)}>> \eps_*\C'_{[0]}
\end{CD}
\end{equation}
is cartesian. For example, all arrows in \eqref{red.sq} are Segal
functors, and the horizontal arrows are $2$-fully faithful. By the
Segal condition, a Segal functor $\gamma$ is $2$-fully faithful if
and only if the square
\begin{equation}\label{ff.1.sq}
\begin{CD}
\C_{[1]} @>{\gamma_{[1]}}>> \C'_{[1]}\\
@V{s^* \times t^*}VV @VV{s^* \times t^*}V\\
\C_{[0]} \times \C_{[0]} @>{\gamma_{[0]} \times \gamma_{[0]}}>>
\C'_{[0]} \times \C'_{[0]}
\end{CD}
\end{equation}
is cartesian, so by \eqref{C.cc}, this corresponds to the usual
notion of a fully faithful functor. In particular, for any functor
$\gamma:I' \to I$, the simplicial replacement $\Delta(\gamma):\Delta
I' \to \Delta I$ is a $2$-functor, and it is $2$-fully faithful iff
$\gamma$ is fully faithful. A {\em natural transformation} between
$2$-functors $\gamma_0,\gamma_1:\C \to \C'$ is a $2$-functor
$\wt{\gamma}:\C \times_{\Delta} \Delta[1] \to \C'$ equipped with
isomorphisms $\wt{\gamma} \circ (\id \times \Delta(\eps(l))) \cong
\gamma_l$ for $l=0,1$. A natural transformation is {\em invertible}
if it extends to $e(\{0,1\}) \supset [1]$, where for any set $S$, we
let $e(S)$ be the category with objects $s \in S$, and exactly one
map between any two objects. A $2$-functor is a {\em
  $2$-equivalence} if it is invertible up to an invertible natural
transformation. A $2$-equivalence is $2$-fully faithful. For the
converse, defining $2$-essentially surjective $2$-functors is
somewhat awkward, but at least, a $2$-fully faithful $2$-functor
$\gamma:\C \to \C'$ that is essentially surjective in the usual
sense is a $2$-equivalence (and by \eqref{ff.sq}, it suffices to
require that $\gamma$ is essentially surjective over $[0] \in
\Delta$).

\section{Abelian categories.}\label{ab.sec}

\subsection{Complexes of length $2$.}

Assume given an abelian category $\A$, and let $C_{[0,1]}(\A)$ be
the category of chain complexes in $\A$ of length $2$ sitting in
homological degrees $0$ and $1$. Any such complex $M_\idot$ defines
a four-term exact sequence
\begin{equation}\label{4.t}
\begin{CD}
0 @>>> H_1(M_\idot) @>>> M_1 @>{d}>> M_0 @>>> H_0(M_\idot) @>>> 0,
\end{CD}
\end{equation}
where $d$ is the differential, and $H_\idot(-)$ stands for homology
objects. Equivalently, $C_{[0,1]}(\A) \cong [1]^o\A$, with a complex
$M_\idot$ corresponding to the functor $[1]^o \to \A$ that sends $l
\in \{0,1\} = [1]$ to $M_l$, and the order relation $0 \leq 1$ to
the differential $d$. The category $C_{[0,1]}(\A)$ has finite limits
and colimits, and we have the following elementary observation.

\begin{lemma}\label{triv.le}
Assume given a map $f_\idot:M_\idot \to M'_\idot$ in $C_{[0,1]}(\A)$
corresponding to a commutative square
\begin{equation}\label{C.A.sq}
\begin{CD}
M_0 @>{f_0}>> M'_0\\
@A{d}AA @AA{d'}A\\
M_1 @>{f_1}>> M'_1
\end{CD}
\end{equation}
in $\A$ such that $f_1$ is surjective. Then $f_\idot$ is a
quasiisomophism if and only if \thetag{i} \eqref{C.A.sq} is
cocartesian and \thetag{ii} $H_1(f_\idot):H_1(M_\idot) \to
H_1(M'_\idot)$ is injective, and if these equivalent conditions
hold, then $f_0$ is surjective.
\end{lemma}

\proof{} Clear. \endproof

For another description of $C_{[0,1]}(\A)$, embed it as a full
subcategory into the category $C_{\geq 0}(\A)$ of all chain
complexes in $\A$ concentrated in non-negative homological degrees,
and the recall that the normalized chain complex functor provides a
{\em Dold-Kan equivalence}
\begin{equation}\label{dk.eq}
C_\idot:\Delta^o\A \cong C_{\geq 0}(\A).
\end{equation}
Then a complex $M_\idot \in C_{\geq 0}(\A)$ lies in $C_{[0,1]}(\A)
\subset C_{\geq 0}(\A)$ iff the corresponding simplicial object $M
\in \Delta^o\A$ satisfies the Segal condition. This is trivially the
case when $M$ is constant --- that is, factors through the oppposite
to the tautological projection $\gamma:\Delta \to \ppt$; this
corresponds to complexes concentrated in degree $0$. In fact,
$\gamma^{o*} \cong \eps^o_!$, where $\eps:\ppt \to \Delta$ is the
embedding onto the terminal object $[0] \in \Delta$; by adjunction,
$\lim_{\Delta^o} = \gamma^o_* \cong \eps^{o*}$ is then given by
evaluation at $[0]$, and since $\gamma \circ \eps \cong \id$,
$\gamma^{o*}:\A \to \Delta^o\A \cong C_{\geq 0}(\A)$ is a fully
faithful embedding. In terms of the equivalence \eqref{dk.eq},
$\eps^{o*}$ sends a complex $M_\idot$ to $M_0$. Again by adjunction,
$\colim_{\Delta^o} = \gamma^o_!$ sends a complex $M_\idot \in
C_{\geq 0}(\A)$ to the cokernel $M_0/d(M_1)$ of the differential
$d:M_1 \to M_0$.

We can also consider the right Kan extension $\eps^o_*:\A \to
\Delta^o\A \cong C_{\geq 0}(\A)$. This also factors through
$C_{[0,1]}(\A) \subset C_{\geq 0}(\A)$ --- for any $M \in \A$, the
object $\eps^o_*M \in \Delta^o\A$ satisfies the Segal condition and
corresponds to the length-$2$ complex $\id:M \to M$.

In homological algebra, exact sequences \eqref{4.t} represent by
Yoneda classes in $\Ext^2(H_0(M_\idot),H_1(M_\idot))$; two sequences
represent the same class iff they are related by a zigzag of maps
equal to $\id$ on $H_0$ and $H_1$. These can be composed with
morphisms in $\A$ in the usual way: for any $N \in \A$ and morphism
$f:H_1(M_\idot) \to N$, we can define $M_\idot \circ f \in
C_{[0,1]}(\A)$ as the complex
\begin{equation}\label{ex.compo}
\begin{CD}
(M_1 \oplus N)/H_1(M_\idot) @>{d \oplus 0}>> M_0,
\end{CD}
\end{equation}
where the differential acts via the projection $M_1 \oplus N \to
M_1$, and then we have natural identifications $H_0(M_\idot \circ
f)\cong H_0(M_\idot)$, $H_1(M_\idot \circ f) \cong N$. Dually, for
any map $g:N \to H_0(M_\idot)$, we have the complex $g \circ
M_\idot$ with identifications $H_1(g \circ M_\idot) \cong
H_1(M_\idot)$, $H_0(g \circ M_\idot) \cong N$ obtained by taking
\eqref{ex.compo} in the opposite category $\A^o$.

By definition, a {\em splitting} of a sequence \eqref{4.t} is given
by an object $\wM \in \A$ equipped with maps
\begin{equation}\label{spl.dia}
\begin{CD}
M_1 @>{a}>> \wM @>{b}>> M_0
\end{CD}
\end{equation}
such that $a$ is injective, $b$ is surjective, and $b \circ a =
d$. A maps of splittings is a map of diagrams \eqref{spl.dia} equal
to $\id$ on $M_0$ and $M_1$. Every such map is invertible, so
splittings form a groupoid $\Spl(M_\idot)$. It is non-empty iff
\eqref{4.t} represents $0$ by Yoneda, and in this case,
$\Spl(M_\idot)$ is non-canonically equivalent to the groupoid
$\Ex^1(H_0(M_\idot),H_1(M_\idot))$ of extensions of $H_0(M_\idot)$
by $H_1(M_\idot)$ (the full truth is that $\Spl(M_\idot)$ is a gerb
over this groupoid of extensions, but we will not need this). A map
$f:H_1(M_\idot) \to N$ defines a functor $f \circ -:\Spl(M_\idot)
\to \Spl(M_\idot \circ f)$, $\wt{M} \mapsto (\wt{M} \oplus
N)/H_1(M_\idot)$, and dually, a map $g:N \to H_0(M_\idot)$ defines a
functor $- \circ g:\Spl(M_\idot) \to \Spl(g \circ M_\idot)$. For any
quasiisomorphism \eqref{C.A.sq} and splitting \eqref{spl.dia} of the
complex $M_\idot$, we can define $\wM'$ as the cokernel of the map
$a \oplus f_1:M_1 \to \wt{M} \oplus M_1'$; this provides a splitting
$$
\begin{CD}
M'_1 @>{0 \oplus \id}>> \wM' = (\wM \oplus M_1')/M_1 @>{(f_0 \circ
  b) \oplus (-d')}>> M_0'
\end{CD}
$$
of the complex $M'_\idot$, functorial with respect to
\eqref{spl.dia}, and the resulting functor $\Spl(M_\idot) \to
\Spl(M'_\idot)$ is an equivalence.

Alternatively, say that a morphism $f:M'_\idot \to M_\idot$ in
$C_{[0,1]}(\A)$ is {\em tight} if $f_1$ is invertible and $H_0(f)$
is injective. Then a splitting \eqref{spl.dia} defines a complex
$M'_\idot$ with $M'_1 = M_1$, $M'_1 = \wM$ and the differential $a$,
we have $H_1(M'_\idot) = 0$, and $b$ provides a tight morphism
$b':M'_\idot \to M_\idot$ with bijective $H_0(b')$. Conversely, any
tight morphism $b':M'_\idot \to M_\idot$ with $H_1(M'_\idot) = 0$
and bijective $H_o(b')$ comes from a unique splitting. More
generally, let $C_{[0,1]}(\A) /^t M_\idot \subset C_{[0,1]}(\A) /
M_\idot$ be the full subcategory spanned by tight morphisms
$M'_\idot \to M_\idot$, and for any object $N \in \A$, let $\A /^i N
\subset \A / N$ be the full subcategory spanned by injective maps
$N' \to N$ (if $N$ only has a set of subobjects, then $\A /^i N$ is
this set partially ordered by inclusion). Then for any tight
$f:M'_\idot \to M_\idot$, $H_1(f)$ is automatically injective, and
$H_0$, $H_1$ define functors $h_l:C_{[0,1]}(\A) /^t M_\idot \to \A
/^i H_l(M_\idot)$, $l=0,1$. The functor $h_0$ is a fibration, the
functor $h_1$ is a cofibration, and $h_0 \times h_1$ has fibers
\begin{equation}\label{spl.cof}
(C_{[0,1]}(\A) /^t M_\idot)_{N_0,N_1} \cong \Spl(q_1 \circ M_\idot
  \circ e_0), \qquad N_l \subset M_l,l = 0,1,
\end{equation}
where $e_l:N_l \to M_l$, $l=0,1$ are the embeddings, and $q_l:M_l
\to M_l/N_l$ are the quotient maps.

Finally, for any abelian category $\A'$, an exact functor $F:\A \to
\A'$ defines a functor
\begin{equation}\label{F.spl}
\Spl(M_\idot) \to \Spl(F(M_\idot))
\end{equation}
obtained by applying $F$ to \eqref{spl.dia} termwise. In general,
\eqref{F.spl} is not an equivalence; however, we have the following
simple observation.

\begin{lemma}\label{spl.le}
Assume that $F$ admits a fully faithful exact left-adjoint functor
$F':\A' \to \A$, and $H_0(M_\idot) \cong F'(N)$ for some $N \in
\A'$. Then \eqref{F.spl} is an equivalence of categories.
\end{lemma}

\proof{} Let $e_\idot:F'(F(M_\idot)) \to M_\idot$ be the adjunction
map. Then $e_\idot$ factors through a quasiisomorphism
$e'_\idot:H_1(e_\idot) \circ F'(F(M_\idot)) \to M_\idot$, and an
inverse equivalence to \eqref{F.spl} is given by the composition
$$
\begin{CD}
\Spl(F(M_\idot)) @>{F'}>> \Spl(F'(F(M_\idot))) @>{H_1(e_\idot) \circ
  -}>> \\ @>{H_1(e_\idot) \circ -}>> \Spl(H_1(e_\idot) \circ
F'(F(M_\idot))) @>>> \Spl(M_\idot),
\end{CD}
$$
where the last functor is the equivalence induced by $e'_\idot$.
\endproof

\subsection{Admissible functors.}\label{adm.subs}

By abuse of notation, for any simplicial set $X \in \Delta^o\Sets$,
with its category of simplices $\Delta X$, let $\eps:X([0]) \to
\Delta X$ be the embedding of the fiber $X([0]) = (\Delta X)_{[0]})$
(treated as a discrete category). For any category $\E$, say that a
functor $E:(\Delta X)^o \to \E$ is {\em special} if the map
$E(\langle [n],x \rangle) \to E(\langle [0],s^*x \rangle)$ is
invertible for any $\langle [n],x \rangle \in \Delta
X$. Equivalently, say that a map $f:[m] \to [n]$ in $\Delta$ is {\em
  special} if $f(0)=0$; then $E:(\Delta )^o \to \E$ is special iff
$E(\langle [n],x \rangle) \to E(\langle [m],f^*x \rangle)$ is
invertible for any $\langle [n],x \rangle \in \Delta X$ and special
$f:[m] \to [n]$.

\begin{exa}\label{xi.i}
For any small category $I$, we have a functor $\xi:(\Delta I)^o \to
I$ sending $\langle [n],i \rangle \in \Delta I$ to $i(0) \in
I$. This functor is special.
\end{exa}

\begin{exa}\label{xi.n}
Let $I=[n]$, for some $n \geq 0$, so that $\Delta[n] \cong \Delta /
[n]$ is the category of objects $[m] \in \Delta$ equipped with a map
$f:[m] \to [n]$. Then the special functor $\xi:(\Delta[n])^o \to
[n]$ of Example~\ref{xi.i} sends $\langle [m],f \rangle \in
\Delta[n]$ to $f(0)$, and it has a left-adjoint $\xi_\dg:[n] \to
(\Delta[n])^o$, $l \mapsto \langle [n-l],t \rangle$ such that the
adjunction map $\id \to \xi \circ \xi_\dg$ is an
isomorphism. Therefore for any target category $\E$, the pullback
functor $\xi^*:\Fun([n],\E) \to (\Delta[n])^o\E$ is a fully faithful
embedding onto the subcategory of functors $E:(\Delta[n])^o \to \E$
such that the adjunction map $\xi^*\xi_\dg^*E \to E$ is an
isomorphism. Since $\xi$ is special, and so are the adjunction maps
$\xi_\dg(\xi(\langle [m],f\rangle)) \to \langle [m],f \rangle$ for
all $\langle [m],f \rangle \in \Delta[n]$, this happens iff $E$
itself is special. Therefore $\xi^*$ identifies $\Fun([n],\E)$ with
the full subcategory in $(\Delta[n])^o\E$ spanned by special
functors, or in other words, a special functor $E:(\Delta[n])^o \to
\E$ uniquely and functorialy factors through $\xi$. In particular,
if $n=0$, then $E:\Delta^o \to \E$ is special iff it is constant.
\end{exa}

\begin{remark}
In fact, for any small $I$, the pullback $\xi^*:\Fun(I,\E) \to
(\Delta I)^o\E$ is a fully faithful embedding onto the full
subcategory spanned by special functors; this fact is quite
well-known but requires a proof (see e.g.\ \cite{DHKS}, or
\cite[Lemma 4.2.1.1]{ka.big}). If $I=[n]$, $\xi$ admits an adjoint,
and the claim becomes obvious.
\end{remark}

\begin{defn}\label{adm.def}
For any $X \in \Delta^o\Sets$, a functor $M_\idot:(\Delta X)^o \to
C_{[0,1]}(\A)$ is {\em admissible} if
\begin{enumerate}
\item the adjunction map $M_1 \to \eps^o_*\eps^{o*}M_1$ is an
  isomorphism, and
\item the functors $H_l(M_\idot):(\Delta X)^o \to \A$, $l=0,1$ are
  special.
\end{enumerate}
\end{defn}

The right Kan extension $\eps^o_*$ in
Definition~\ref{adm.def}~\thetag{i} can be computed explicitly by
\eqref{kan.eq}. This shows that $M_\idot$ satisfies
Definition~\ref{adm.def}~\thetag{i} if and only if for any $\langle
[n],x \rangle \in \Delta X$, with the embeddings $\eps_m:[0] = \ppt
\to [n]$, $0 \mapsto m$ for all $m \in [n]$, the map
\begin{equation}\label{e.l}
M_1(\langle [n],x \rangle) \to \bigoplus_{m \in [n]} M_1(\langle
[0],\eps_m^*x \rangle)
\end{equation}
is an isomorphism. In particular, for any $Y \in \Delta^o\Sets$ and
map $f:Y \to X$, $\Delta(f)^*$ sends admissible functors to
admissible functors. Moreover, for any admissible $M_\idot:(\Delta
X)^o \to C_{[0,1]}(\A)$ and $\langle [n],x \rangle \in \Delta$, the
quasiisomorphism $M_\idot(\langle [n],x \rangle) \to M_\idot(\langle
[0],s^*x \rangle)$ is surjective (for $M_1$, this immediately
follows from \eqref{e.l}, and then for $M_0$, this is
Lemma~\ref{triv.le}).

\begin{exa}\label{seg.bi.exa}
Take $X=\ppt$. By Example~\ref{xi.n}, a functor $E:\Delta^o \to \E$
to any category $E$ is special iff it is constant. A functor
$M_\idot:\Delta^o \to C_{[0,1]}(\A)$ admissible in the sense of
Definition~\ref{adm.def} automatically satisfies the Segal
condition, thus corresponds to a bicomplex
\begin{equation}\label{seg.bi}
\begin{CD}
C_1(M_0) @>>> C_0(M_0)\\
@A{b}AA @AAA\\
C_1(M_1) @>{a}>> C_0(M_1)
\end{CD}
\end{equation}
whose rows are given by \eqref{dk.eq}. Then for any
$M_\idot:\Delta^o \to C_{[0,1]}(\A)$ satisfying the Segal condition,
with the correponding bicomplex \eqref{seg.bi},
Definition~\ref{adm.def}~\thetag{i} resp.\ \thetag{ii} means that
$a$ resp.\ $b$ is invertible.
\end{exa}

\begin{prop}\label{seg.prop}
For any integers $n > l > 0$, with the corresponding embedding
\eqref{seg.eq}, the pullback functor $\Delta(b^l_n)^*$ induces an
equivalence between the full subcategories spanned by admissible
functors to $C_{[0,1]}(\A)$.
\end{prop}

\proof{} To simplify notation, let $N([n])_l$ be the source of the
map \eqref{seg.eq}, and denote $(\Delta [n])_l = \Delta N([n])_l$;
explicitly, $(\Delta [n])_l \subset \Delta [n] = \Delta / [n]$ is
the full subcategory of arrows $f:[m] \to [n]$ that factor through
$s([l]) \subset [n]$ or $t([n-l]) \subset [n]$. Let $(\Delta [n])^l
\subset \Delta[n]$ be the full subcategory of $f:[m] \to [n]$ with
$l \in f([m])$, and let $(\Delta[n])_{l\dm} = (\Delta[n])_l \cup
(\Delta[n])^l \subset \Delta[n]$. Then the full embedding
$\Delta(b^l_n):(\Delta[n])_l \to \Delta[n]$ factors as
\begin{equation}\label{ab.eq}
\begin{CD}
(\Delta[n])_l @>{\alpha}>> (\Delta[n])_{l\dm} @>{\beta}>>
  \Delta[n],
\end{CD}
\end{equation}
where both functors are full embeddings. For any object $\langle
[n],f \rangle \in (\Delta[n])_{l\dm}$, the object $\langle [0],s^*f
\rangle \in \Delta[n]$ lies in $\Delta([n])_{l\dm} \subset
\Delta[n]$, so that Definition~\ref{adm.def}~\thetag{ii} makes sense
for functors $(\Delta[n])_{l\dm}^o \to C_{[0,1]}(\A)$. More
generally, $\langle [0],\eps_m^*f \rangle$ is in $(\Delta[n])_{l\dm}$
--- in fact, already in $(\Delta[n])_l$ --- for any $m \in [n]$, so
Definition~\ref{adm.def}~\thetag{i} in the form \eqref{e.l} also
makes sense, and we can speak about admissible functors
$(\Delta[n])_{l\dm}^o \to C_{[0,1]}(\A)$. For any admissible
$M_\idot:(\Delta [n])^o \to C_{[0,1]}(\A)$, $\beta^{o*}M_\idot$ is
admissible, and for any admissible $M_\idot:(\Delta[n])_{l\dm}^o \to
C_{[0,1]}(\A)$, $\alpha^{o*}M_\idot$ is also admissible.

\begin{lemma}\label{adm.le}
For any admissible functor $M_\idot:(\Delta[n])^o_l \to \A$, the
functor $M'_\idot=\alpha^o_*M_\idot$ is admissible. For any
admissible $M_\idot:(\Delta[n])_{l\dm}^o \to \A$, the functor
$M'_\idot=\beta^o_!M_\idot$ is admissible as well.
\end{lemma}

Assuming this for the moment, we can finish the proof. Indeed, since
$\alpha$ and $\beta$ in \eqref{ab.eq} are fully faithful, we have
$M_\idot \cong \alpha^{o*}\alpha^o_*M_\idot \cong (\beta \circ
\alpha)^{o*}\beta^o_!\alpha^o_*M_\idot$ for any
$M_\idot:((\Delta[n])_l \times [1])^o \to \A$. On the other hand, we
have adjunction maps
\begin{equation}\label{ad.eq}
M'_\idot \to \alpha^o_*\alpha^{o*}M_\idot', \qquad
\beta^o_!\beta^{o*}M_\idot \to M_\idot
\end{equation}
for any $M_\idot:(\Delta[n])^o \to C_{[0,1]}(\A)$,
$M'_\idot:(\Delta[n])_{l\dm}^o \to C_{[0,1]}(\A)$, and once we know
that all the functors in \eqref{ad.eq} are admissible, \eqref{e.l}
and Lemma~\ref{triv.le} immediately imply that both maps are
isomorphisms. Therefore $\beta^o_!\alpha^o_*$ provides an inverse
equivalence to $\Delta(b_{l,n})^{o*} \cong \alpha^{o*}\beta^{o*}$.

\proof[Proof of Lemma~\ref{adm.le}.] In the first claim, since
$\alpha$ is an equivalence over $[0] \in \Delta$,
Definition~\ref{adm.def}~\thetag{i} immediately follows from
\eqref{e.l}. For Definition~\ref{adm.def}~\thetag{ii}, we need to
check that the functor $\alpha^o_*M_\idot:(\Delta[n])_{l\dm}^o \to
C_{[0,1]}(\A)$ sends maps $\langle [0],s^*f \rangle \to \langle
[m],f \rangle$ to quasiisomorphisms. Since $\alpha$ is fully
faithful and $M_\idot$ is admissible, the claim holds if $\langle
[m],f\rangle \in (\Delta[n])_l$, so assume $\langle [m],f \rangle
\in (\Delta[n])^l$. Then $f^{-1}(\{l\}) \subset [m] = \{0,\dots,m\}$
is an interval $\{p,\dots,q\}$, we have an embedding $\V \to
(\Delta[n])_l / \langle [m],f \rangle$ represented by the diagram
$$
\begin{CD}
\langle [q],s^*f \rangle @<<< \langle [q-p],s^*t^*f \rangle @>>>
\langle [m-p],t^*f \rangle,
\end{CD}
$$
and this embedding admits a left-adjoint, so it is cofinal.
We then have a commutative diagram
$$
\begin{CD}
\alpha^o_*M_\idot(\langle [m],f \rangle) @>{b'}>> M_\idot(\langle
   [q],s^*f \rangle) @>>> M_\idot(\langle [0],s^*f \rangle\\
@VVV @VVV\\
M_\idot(\langle [m-p],t^*f \rangle) @>{b}>> M_\idot(\langle
[q-p],s^*t^*f \rangle),
\end{CD}
$$
since $M_\idot$ is admissible, the maps $a$ and $b$ are surjective
quasiisomorphisms, and the square is cartesian by \eqref{kan.eq}, so
that $b'$ is also a quasiisomorphism.

In the second claim, analogously, it suffices to consider $\langle
[m],f \rangle \in \Delta[n]$ that is not in
$(\Delta[n])_{l\dm}$. Explicitly, the comma-fiber $\langle [m],f
\rangle \setminus (\Delta[n])_{l\dm}$ is the category of
factorizations
\begin{equation}\label{m.f.facto}
\begin{CD}
[m] @>{e}>> [m_+] @>{f_+}>> [n]
\end{CD}
\end{equation}
of the map $f:[m] \to [n]$ such that $l \in f_+([m_+])$. Say that a
factorization \eqref{m.f.facto} is {\em tight} if $e$ is injective
and $[m_+] = f_+^{-1}(l) \cup e([m])$, and note that tight
factorizations form a full subcategory in $\langle [m],f \rangle
\setminus (\Delta[n])_{l\dm}$ equivalent to $\Delta$, with the
equivalence $\langle [m_+],f_+ \rangle \mapsto f^{-1}(l)$. The
corresponding full embedding $\nu':\Delta \subset \langle [m],f
\rangle \setminus (\Delta[n])_{l\dm}$ admits a right-adjoint, so it
is final. Thus if we denote $\nu = \sigma(\langle [n],f \rangle) \circ
\nu':\Delta \to \Delta[n]$, then \eqref{kan.eq} provides an
isomorphism $\beta^o_!M_\idot(\langle [m],f \rangle) \cong
\colim_{\Delta^o}\nu^{o*}M_\idot$. This colimit can be computed as
the cokernel of the differential $C_1(\nu^{o*}M_\idot) \to
C_0(\nu^{o*}M_\idot)$ in the normalized chain complex
\eqref{dk.eq}. Moreover, since $M_\idot$ is admissible, while $\nu$
sends squares \eqref{seg.sq} resp.\ special maps to squares
\eqref{seg.sq} resp.\ special maps, $\nu^{o*}M_\idot$ satisfies the
Segal condition, and the map $b$ in the corresponding bicomplex
\eqref{seg.bi} is an isomorphism. Therefore if we identify
$C_0(\nu^{o*}M_\idot) \cong M_\idot(\nu([0]))$, then \eqref{seg.bi}
gives rise to a commutative diagram
\begin{equation}\label{be.dia}
\begin{CD}
C_1(\nu^{o*}(M_1)) @>{a}>> M_1(\nu([0])) @>>> M_0(\nu([0]))\\
@VVV @VV{q_1}V @VV{q_0}V\\
0 @>>> \beta^o_!M_1(\langle [m],f \rangle) @>>> \beta^o_!M_0(\langle
[m],f \rangle)
\end{CD}
\end{equation}
with cocartesian squares. The isomorphism \eqref{e.l} for $M_1$ then
reads as
$$
M_1(\nu([0])) \cong M_1(\langle [0],\eps_l \rangle) \oplus \bigoplus_{p
  \in [m]} M_1(\langle [0],\eps_{f(p)} \rangle),
$$
and in terms of this isomorphism, $a$ in \eqref{be.dia} is the
embedding onto the first summand, so we obtain \eqref{e.l}, hence
Definition~\ref{adm.def}~\thetag{i} for $\beta^o_!M_1(\langle [m],f
\rangle)$. To finish the proof, it remains to check that $q_\idot$
in \eqref{be.dia} is a quasiisomorphism, and by Lemma~\ref{triv.le},
this amounts to checking that $H_1(q_\idot)$ is injective. However,
since $M_\idot$ is admissible, the map $g_\idot:M_\idot(\nu([0]))
\to M_\idot(\langle [0],s^*f \rangle)$ is a surjective
quasiisomorphism by assumption, so $H_1(g_\idot)$ is injective by
Lemma~\ref{triv.le}, and $H_1(g_\idot)$ factors through
$H_1(q_\idot)$.
\endproof

\begin{remark}
Our proof of Proposition~\ref{seg.prop} by using the decomposition
\eqref{ab.eq} is also quite standard; see e.g. \cite[Subsection
  8.3]{ka.bo}.
\end{remark}

\subsection{Constructing the $2$-category.}\label{2.co.subs}

We can now construct our $2$-category of extensions
$C^{(2)}_{[0,1]}(\A)$. Consider the twisted arrow category $\tw(\Delta)$,
with the fibration $\tau:\tw(\Delta) \to \Delta^o$ of
Example~\ref{tw.exa}, and denote by
\begin{equation}\label{wt.c}
  \wt{\C}_{[0,1]}(\A) \subset \Fun(\tw(\Delta)^o|\Delta,C_{[0,1]}(\A))
\end{equation}
the full subcategory of functors $\tw(\Delta)^o_{[n]} \cong (\Delta
/ [n])^o \cong (\Delta[n])^o \to C_{[0,1]}(\A)$ that are admisible
in the sense of Definition~\ref{adm.def}. Then
Proposition~\ref{seg.prop} immediately implies the following.

\begin{corr}\label{seg.corr}
The category \eqref{wt.c} is a Segal category.
\end{corr}

\proof{} The relative functor category in \eqref{wt.c} is by
definition fibered over $\Delta$, with fibers
$(\Delta[n])^oC_{[0,1]}(\A)$; the transition functor for a map
$f:[m] \to [n]$ is given by the pullback with respect to
$\Delta(f)^o:(\Delta[m])^o \to (\Delta[n])^o$. Since
$\Delta(f)^{o*}$ sends admissible functors to admissible functors,
the induced projection $\wt{C}_{[0,1]}(\A) \to \Delta$ is also a
fibration, and the Segal property immediately follows from
\eqref{IX.eq} and Proposition~\ref{seg.prop}.
\endproof

\begin{defn}
The {\em extension $2$-category} $C^{(2)}_{[0,1]}(\A)$ of the
abelian category $\A$ is the reduction $\wt{\C}_{[0,1]}(\A)^{\red}$
of the Segal category \eqref{wt.c} given by \eqref{red.sq}.
\end{defn}

By definition, objects of the $2$-category $C^{(2)}_{[0,1]}(\A)$ are
isomorphism classes of admissible functors $\Delta^o \to
C_{[0,1]}(\A)$, and by Example~\ref{seg.bi.exa}, this is the same
thing as isomorphism classes of complexes $M_\idot \in
C_{[0,1]}(\A)$. To describe the categories of morphisms, assume
given two objects $M_\idot,N_\idot \in C_{[0,1]}(\A)$, and let
$C^{(2)} = C^{(2)}_{[0,1]}(\A)$ to simplify notation. Then by
definition, $C^{(2)}(M_\idot,N_\idot)$ is the category of admissible
functors $F_\idot:(\Delta[1])^o \to C_{[0,1]}(\A)$ equipped with
isomorphisms
\begin{equation}\label{F.MN}
(\Delta(\eps_0)^o \times \Delta(\eps_1)^o)^*F_\idot \cong M^\Delta_\idot
  \times N^\Delta_\idot,
\end{equation}
where $M^\Delta_\idot,N^\Delta_\idot:\Delta^o \to C_{[0,1]}(\A)$ are
the admissible functors corresponding to $M_\idot,N_\idot \in
C_{[0,A]}$, and $\eps_l:\ppt \to [1]$, $l=0,1$ are the embeddings
onto $l$. By Definition~\ref{adm.def}~\thetag{ii}, such an
admissible functor $F_\idot$ gives rise to special functors
$H_l(F_\idot):(\Delta[1])^o \to \A$, $l=0,1$, and by
Example~\ref{xi.n}, we have $H_l(F_\idot) \cong
\xi^*\xi_*H_l(F_\idot)$, where $\xi_*H_l(F_\idot) \in \Fun([1],\A)$,
$l = 0,1$ are functors $[1] \to \A$. Moreover, the isomorphism
\eqref{F.MN} provides isomorphisms $\xi_*H_l(F_\idot)(0) \cong
H_l(M_\idot)$, $\xi_*H_l(F_\idot)(1) \cong H_l(N_\idot)$, so the
only remaining data are maps $f_l:H_l(M_\idot) \to H_l(N_\idot)$,
$l=0,1$. This defines a functor
\begin{equation}\label{C.MN}
C^{(2)}(M_\idot,N_\idot) \to \Hom(H_0(M_\idot),H_0(N_\idot)) \times
\Hom(H_1(M_\idot),H_1(N_\idot))
\end{equation}
whose target is understood as a discrete category (that is,
$C^{(2)}(M_\idot,N_\idot)$ decomposes into a disjoint union of its
fibers over all $f_0 \times f_1$).

\begin{prop}\label{cat.prop}
For any pair of maps $f_l:H_l(M_\idot) \to H_l(N_\idot)$, $l=0,1$,
the fiber of the projection \eqref{C.MN} over $f_0 \times f_1$ is
given by
\begin{equation}\label{spl.3}
C^{(2)}(M_\idot,N_\idot)_{f_0 \times f_1} \cong \Spl((f_1 \oplus
\id) \circ (M_\idot \oplus N_\idot) \circ (\id \oplus (-f_0))).
\end{equation}
\end{prop}

\proof{} Functors $A:[1] \to \A$ are the same thing as arrows $f:A_0
\to A_1$ in $\A$. For any such $A$, the adjunction map
$e^\Delta:\xi^*A \to \Delta(\eps_0)^o_*\xi^*A_0 \oplus
\Delta(\eps_1)^o_*\xi^*A_1$ is injective --- indeed, since the
embedding $\eps:\{0,1\} \to \Delta[1]$ factors through
$\Delta(\eps_0) \copr \Delta(\eps_1):\Delta \copr \Delta \to
\Delta[1]$, the adjunction map $e':\xi^*A \to \eps_*\eps^*\xi^*A$
factors through $e^\Delta$, and $e'$ is injective by
\eqref{e.l}. Therefore we have a functorial short exact sequence
\begin{equation}\label{nu.1}
0 \longrightarrow \xi^*A \overset{e^\Delta}\longrightarrow
\Delta(\eps_0)^o_*\xi^*A_0 \oplus \Delta(\eps_1)^o_*\xi^*A_1
\overset{q^\Delta}\longrightarrow K(A) \longrightarrow 0,
\end{equation}
where $K:\Fun([1],\A) \to (\Delta[1])^o\A$ is a certain functor
whose precise form is not important. What is important is that the
Kan extension $\xi_*$ is exact by Example~\ref{xi.n}, and if we
apply it to \eqref{nu.1}, we obtain a functorial short exact
sequence
\begin{equation}\label{nu.2}
\begin{CD}
0 @>>> A @>{e}>> \eps_{0*}A_0 \oplus \eps_{1*}A_1 @>{q}>> \eps_{0*}A_1
@>>> 0,
\end{CD}
\end{equation}
where explicitly, by \eqref{kan.eq}, $\eps_{0*}A_0$ corresponds to
$A_0 \to 0$, $\eps_{1*}A_1$ corresponds to $\id:A_1 \to A_1$, and
the quotient $\xi_*K(A)$ then corresponds to $A_1 \to 0$ and can be
identified with $\eps_{0*}A_1$. The arrow $e$ in \eqref{nu.2} is
$\id \oplus f$ resp.\ $\id$ at $0 \in [1]$ resp.\ $1 \in [1]$, and
the arrow $q$ is $(-f) \oplus \id$ resp.\ $0$.

In particular, for any admissible $F_\idot:(\Delta[1])^o \to
C_{[0,1]}(\A)$, the functors $\xi_*H_l(F_\idot):[1] \to \A$, $l=0,1$
give rise to the sequences \eqref{nu.1}, \eqref{nu.2}; let us
denote the corresponding maps $e^\Delta$, $q^\Delta$, $e$, $q$ by
$e^\Delta_l$, $q^\Delta_l$, $e_l$, $q_l$, $l=0,1$. Then by
adjunction, giving a map \eqref{F.MN} is equivalent to giving a map
\begin{equation}\label{g.MN}
g_\idot:F_\idot \to \Delta(\eps_0)^o_*M_\idot^\Delta \oplus
\Delta(\eps_1)^o_*N_\idot^\Delta,
\end{equation}
and \eqref{F.MN} is an isomorphism iff the corresponding map
\eqref{g.MN} is tight, and $H_l(g_\idot) = e_l$, $l=0,1$. Then
\eqref{spl.cof} provides an equivalence
\begin{equation}\label{spl.1}
  C^{(2)}(M_\idot,N_\idot)_{f_0 \times f_1} \cong \Spl(q^\Delta_1
  \circ (\Delta(\eps_0)^o_*M_\idot^\Delta \oplus
  \Delta(\eps_1)^o_*N_\idot^\Delta) \circ e^\Delta_0).
\end{equation}
Moreover, the adjoint pair $\xi_*$, $\xi^*$ lies within the scope of
Lemma~\ref{spl.le}, and this reduces \eqref{spl.1} to an equivalence
\begin{equation}\label{spl.2}
C^{(2)}(M_\idot,N_\idot)_{f_0 \times f_1} \cong \Spl(q_1 \circ
(\eps_{0*}M_\idot \oplus \eps_{1*}N_\idot) \circ e_0).
\end{equation}
However, while $q_1 \circ (\eps_{0*}M_\idot \oplus
\eps_{1*}N_\idot) \circ e_0$ in \eqref{spl.2} is a length-$2$
complex in the category $\Fun([1],\A)$, its degree-$1$ homology is
of the form $\eps_{0*}H_1(N_\idot)$ by virtue of \eqref{nu.2}. We
can then dualize Lemma~\ref{spl.le} by passing to the opposite
abelian categories, and apply it to the adjoint pair $\eps_{0*}$,
$\eps_0^*$. Since $\eps_0^*\eps_{0*}\cong\eps_0^*\eps_{1*}\cong\id$,
this gives \eqref{spl.3} on the nose.
\endproof

Explicitly, if the four-term sequences \eqref{4.t} corresponding to
our complexes $M_\idot,N_\idot \in C_{[0,1]}(\A)$ represent classes
$\mu \in \Ext^2(H_0(M_\idot),H_1(M_\idot))$ and $\nu \in
\Ext^2(H_0(N_\idot),H_1(N_\idot))$, then $(f_1 \oplus \id) \circ
(M_\idot \oplus N_\idot) \circ (\id \oplus (-f_0))$ in \eqref{spl.3}
represents the class $\nu \circ f_1 - f_0 \circ \mu \in
\Ext^2(H_0(M_\idot),H_1(N_\idot))$. Proposition~\ref{cat.prop} then
says that $C^{(2)}(M_\idot,N_\idot)_{f_0 \times f_1}$ is non-empty
iff $\nu \circ f_1 = f_0 \circ \mu$, and if it is non-empty, then it
is equivalent to the groupoid $\Ex^1(H_0(M_\idot),H_1(N_\idot))$ (or
more precisely, is a gerb over this groupoid). If the abelian
category $\A$ admits a derived category $\D(\A)$, so that
$\Ext^\hdot(-,-)$ in $\A$ are sets and not proper classes, then the
$2$-category $\C^{(2)}_{[0,1]}(\A)$ is bounded, and its truncation
is the full subcategory $\D_{[0,1]}(\A) \subset \D(\A)$ spanned by
objects concentrated in homological degrees $0$ and $1$.

\section{Derived categories.}\label{der.sec}

\subsection{Generalities on enhancements.}

Let $\Pos$ be the category of partially ordered sets. In the
enhancement formalism of \cite{ka.big}, an {\em enhanced category}
is a Grothendieck fibration $\C \to \Pos$ satisfying a number of
axioms. The precise shape of these axioms can be found in
\cite[Subsection 7.2.1]{ka.big} or \cite[Subsection 3.2]{ka.small};
however, it is irrelevant for our present purposes. An {\em enhanced
  functor} between enhanced categories is a functor cartesian over
$\Pos$. The {\em $1$-truncation} of an enhanced category $\C \to
\Pos$ is the fiber $\C_{\ppt}$ over the one-point set $\ppt \in
\Pos$; the whole $\C$ is understood as an enhancement for its
$1$-truncation. Any category $\E$ carries a trivial enhancement
$\Unf(\E) = \Fun(\Pos^\hdot|\Pos,\E)$, where $\Pos^\hdot \to \Pos$
is the tautological cofibration with fibers $\Pos^\hdot_J \cong
J^o$. For any enhanced category $\C$, we have a canonical enhanced
functor $\unf(\C):\C \to \Unf(\C_\ppt)$, and for any enhanced
functor $\gamma:\C \to \Unf(\E)$ for some $\E$, we have a canonical
isomorphism $\gamma \cong \Unf(\gamma_\ppt) \circ \unf(\C)$. We say
that $\C$ is {\em $1$-truncated} if $\unf(\C)$ is an equivalence;
this happens iff $\C \cong \Unf(\E)$ for some $\E$. We denote
$\ppt^h = \Unf(\ppt)$, and this is just $\Pos$ equipped with the
trivial fibration $\id:\Pos \to \Pos$. An enhanced functor is {\em
  small} if it is small as a functor, and an enhanced category $\C$
is {\em small} if so is the structural fibration $\C \to \ppt^h
\cong \Pos$. For any enhanced category $\C$, {\em enhanced objects}
resp.\ {\em morphisms} in $\C$ are objects $c \in \C_\ppt$
resp.\ morphisms in $\C_\ppt$, and any such $c$ defines a unique
enhanced functor $\eps^h(c):\ppt^h \to \C$. An enhanced functor $\C'
\to \C$ is fully faithful if it is fully faithful in the usual
sense, and enhanced full subcategories $\C' \subset \C$ correspond
bijectively to full subcategories $\C'_\ppt \subset \C_\ppt$ (that
is, any full subcategory $\C'_\ppt \subset \C_\ppt$ inherits a
canonical enhancement).

While exact axioms characterizing enhanced categories are not needed
for our applications, we do need the following two fundamental
facts.
\begin{enumerate}
\item For any enhanced category $\E$ and small enhanced category
  $\E$, there exists an {\em enhanced functor category}
  $\fFun^h(\C,\E)$ satisfying the usual universal property
  (\cite[Corollary 7.3.3.5]{ka.big} or \cite[Corollary
  3.2]{ka.small}). If $\C$ is $1$-truncated --- that is,
  $\C=\Unf(I)$ for some small category $I$ --- we will simplify
  notation by writing $\fFun^h(I,\E) = \fFun^h(\C,\E)$. If $I$ is a
  partially ordered set, then $\fFun^h(I,\E) \cong i^*\E$, where
  $i:\Pos \to \Pos$ is given by $i(J) = J \times I^o$.
\item For any small enhanced categories $\C$, $\C_0$, $\C_1$
  equipped with enhanced functors $\C_0,\C_1 \to \C$, there exists a
  small enhanced category $\C_0 \times^h_{\C} \C_1$ equipped with a
  functor $\gamma:\C_0 \times^h_{\C} \C_1 \to \C_0 \times_{\C} \C_1$
  cartesian over $\Pos$ that is an {\em epivalence} --- that is,
  conservative, essentially surjective and full. Moreover, for any
  such $\gamma:\C_0 \times^h_{\C} \C_1 \to \C_0 \times_{\C} \C_1$,
  and any other small enhanced category $\C'$ equipped with a
  functor $\gamma':\C' \to \C_0 \times_{\C} \C_1$ cartesian over
  $\Pos$, there exists an enhanced functor $\wt{\gamma}:\C' \to \C_0
  \times^h_{\C} \C_1$ and an isomorphism $\alpha:\gamma' \cong
  \gamma \circ \wt{\gamma}$, and the pair $\langle
  \wt{\gamma},\alpha \rangle$ is unique up to an isomorphism. These
  results are \cite[Lemma 7.3.3.7]{ka.big} or \cite[Corollary 3.23
    \thetag{i}]{ka.small} (existence), and \cite[Corollary
    7.3.3.6]{ka.big} or \cite[Corollary 3.23 \thetag{ii}]{ka.small}
  (existence and functoriality).
\end{enumerate}
The enhanced category $\C_0 \times^h_{\C} \C_1$ is called the {\em
  semicartesian product} of $\C_0$ and $\C_1$ over $\C$, and the
universal property characterizes it uniquely up to a unique
equivalence. In particular, if $\C$ is $1$-truncated, then the usual
product $\C_0 \times_{\C} \C_1$ is an enhanced category, and we then
have $\C_0 \times^h_{\C} \C_1 \cong \C_0 \times_{\C} \C_1$. The
smallness assumption can be slightly relaxed; for example, it
suffices to require that $\C_1 \to \C$ is small, and $\C_0$ is small
or discrete (for the latter, see \cite[Remark 7.3.3.8]{ka.big}).

Using \thetag{i} and \thetag{ii} above as a black box, one develops
the theory of enhanced categories in a way that is largely parallel
to the usual category theory. In particular, one can distinguish a
class of enhanced functors $\pi:\C \to \E$ between small enhanced
categories that are {\em enhanced fibrations} or {\em enhanced
  cofibrations}. An enhanced fibration $\C \to \E$ has enhanced
fibers $\C_e$ for all enhanced objects $e \in \E_\ppt$, and enhanced
{\em transition functors} $f^*:\C_{e'} \to \C_e$ for all enhanced
morphisms $f:e \to e'$, and dually for enhanced
cofibrations. Moreover, one has an enhanced version of the relative
functor category construction: for any enhanced cofibration
$\E^\hdot \to \E$ between small enhanced categories, and any
enhanced category $\C$, one construct the {\em relative enhanced
  functor category} $\fFun^h(\E^\hdot|\E,\C)$, with the appropriate
universal property. It comes equipped with an enhanced fibration
$\fFun^h(\E^\hdot|\E,\C) \to \E$ with enhanced fibers
$\fFun^h(\E^\hdot|\E,\C)_e \cong \fFun^h(\E^\hdot_e,\C)$, $e \in
\E_\ppt$. In particular, for any cofibration $\E^\hdot \to \E$
between small categories, $\Unf(\E^\hdot) \to \Unf(\E)$ is an
enhanced cofibration; we will simplify notation by writing
$\fFun^h(\E^\hdot|\E,C) = \fFun^h(\Unf(\E^\hdot)|\Unf(\E),\C)$.

For any abelian category $\A$, the chain-homotopy category $\Ho(\A)$
admits a natural enhancement $\Hho(\A)$ obtained as an appropriate
localization of $\Unf(C_\idot(\A))$, see \cite[Lemma
  7.2.4.2]{ka.big}, \cite[Proposition 3.17]{ka.small}. In principle,
provided that $\A$ is nice enough to admit a derived category
$\D(\A)$, this derived category also has a natural
enhancement. However, proving this requires the full theory of
stable enhanced categories and Verdier localization, and this is is
not available at the moment. So for simplicity, let us assume that
$\A$ has enough injectives. In this case $\D(\A)$ definitely exists,
and can be realized as the full subcategory $\D(\A) \subset \Ho(\A)$
spanned by $h$-injective complexes; the quotient functor $q:\Ho(\A)
\to \D(\A)$ is left-adjoint to the embedding $\D(\A) \to
\Ho(\A)$. Then being a full subcategory, $\D(\A) \subset \Ho(\A)$
inherits an enhancement $\D(\A)^h \subset \Hho(\A)$, and the
embedding $\D(\A)^h \subset \Hho(\A)$ admits a left-adjoint enhanced
functor $q^h:\Hho(\A) \to \D(\A)^h$ that enhances the quotient
functor $q$. Composing $q^h$ with the projection $\Unf(C_\idot(\A))
\to \Hho(\A)$ gives an enhanced functor $h:\Unf(C_\idot(\A)) \to
\D(\A)^h$. For any $J \in \Pos$, we have $\Unf(C_\idot(\A))_J \cong
J^oC_\idot(\A) \cong C_\idot(J^o\A)$ and $\D(\A)^h_J \cong
\D(J^o\A)$, while $h:C_\idot(J^o\A) \to \D(J^o\A)$ is the usual
localization functor; overall, $\D(\A)^h$ is obtained by localizing
$\Unf(C_\idot(\A))$ with respect to the class of maps consisting of
quasiisomorphisms in all the fibers $C_\idot(J^o\A) \cong
\Unf(C_\idot(\A))_J \subset \Unf(C_\idot(\A))$ of the fibration
$\Unf(C_\idot(\A)) \to \Pos$. We can further restrict our attention
to the full subcategory $\D_{[0,1]}(\A) \subset \D(\A)$ spanned by
complexes concentrated in homological degrees $0$ and $1$, and we
obtain an enhancement $\D_{[0,1]}(\A)^h$ and an enhanced functor
$h:\Unf(C_{[0,1]}(\A)) \to \D_{[0,1]}(\A)^h$.

\subsection{Enhanced $2$-categories.}

To construct en enhanced version of the machinery of $2$-categories,
one can use exactly the same packaging as in
Subsection~\ref{2cat.subs}. One defines an {\em enhanced Segal
  category} as an enhanced fibration $\C \to \Unf(\Delta)$, with
enhanced fibers $\C_{[n]}$, $[n] \in \Delta$ such that for any
square \eqref{seg.sq}, the corresponding functor
\begin{equation}\label{enh.seg.eq}
\C_{[n]} \to \C_{[l]} \times_{\C_{[0]}} \C_{[n-l]}
\end{equation}
is an epivalence (equivalently, $\C_{[n]} \cong \C_{[l]}
\times^h_{\C_{[0]}} \C_{[n-1]}$). For any enhanced category $\E$, a
trivial example is provided by an enhanced version of \eqref{eps.eq}
that reads as
\begin{equation}\label{eps.h.eq}
\eps_*^h\E = \fFun^h([0] \setminus \Delta|\Delta,\E),
\end{equation}
where the right-hand side is taken as the definition of the
left-hand side. An {\em enhanced $2$-category} is an enhanced Segal
category $\C$ with discrete $\C_{[0]}$. As in
Subsection~\ref{2cat.subs}, for any enhanced $2$-category $\C$, the
transition enhanced functor $s^* \times t^*$ induces a decomposition
\eqref{C.cc}, where $\C(c,c')$ are enhanced categories, small if
$\C$ is bounded. Moreover, since $\C_{[0]}$ is discrete, and in
particular $1$-truncated, the product in \eqref{enh.seg.eq}
coincides with the semicartesian product, so that the truncation
$\C_\ppt$ of the enhanced $2$-category $\C$ is a $2$-category in the
unenhanced sense, with the same objects, and morphism categories
$\C_\ppt(c,c') = \C(c,c')_\ppt$.

Now, as it happens, Example~\ref{repl.exa} also has an enhanced
counterpart. Namely, for any enhanced category $\E$, we can consider
the enhanced category $\Delta^{h\Tot}\E =
\fFun^h(\Delta^\hdot|\Delta,\E)$. Then the embedding $[0] \setminus
\Delta \to \Delta^\hdot$ induces an enhanced functor
$\gamma:\Delta^{h\Tot}\E \to \eps^h_*\E$, one checks that this
enhanced functor is small, and we can then define an enhanced
category $\Delta^h \E$ by the semicartesian square
\begin{equation}\label{red.h.sq}
\begin{CD}
\Delta^h\E @>>> \Delta^{h\Tot}\E\\
@VVV @VV{\gamma}V\\
\eps^h_*\Unf(\pi_0(\E_{\ppt\Iso})) @>>> \eps^h_*\E,
\end{CD}
\end{equation}
where as in the enhanced case, the bottom arrow is obtained by
choosing an enhanced object $e \in \E_\ppt$ in each isomorphism
class. By \cite[Proposition 7.5.6.2]{ka.big}, $\Delta^{h\Tot}\E$ is
an enhanced Segal category, and $\Delta^h\E$ is an enhanced
$2$-category, so that its truncation $(\Delta^h\E)_\ppt$ is a
$2$-category in the usual sense.

\begin{prop}
For any abelian category $\A$ with enough injectives, with the
corresponding enhanced category $\E=\D_{[0,1]}(\A)^h$, the enhanced
$2$-category $\Delta^h\E$ is $1$-truncated and naturally
$2$-equivalent to the $2$-category $\C_{[0,1]}^{(2)}(\A)$ of
Corollary~\ref{seg.corr}.
\end{prop}

\proof{} As in Subsection~\ref{2.co.subs}, consider the twisted
arrow category $\tw(\Delta)$, with the fibration $\tau:\tw(\Delta)
\to \Delta^o$ of Example~\ref{tw.exa}. The opposite functor
$\tau^o:\tw(\Delta)^o \to \Delta$ is then a cofibration with fibers
$\tw(\Delta)^o_{[n]} \cong (\Delta[n])^o$, and the functors
$\xi:(\Delta[n])^o \to [n]$, $[n] \in \Delta$ of Example~\ref{xi.i}
together provide a functor $\xi_\idot:\tw(\Delta)^o \to
\Delta^\hdot$ cocartesian over $\Delta$. For any enhanced category
$\E$, we then have an enhanced functor
\begin{equation}\label{xi.do}
\xi_\idot^*:\Delta^{h\Tot}\E = \fFun^h(\Delta^\hdot)|\Delta,\E) \to
\fFun^h(\tw(\Delta)^o)|\Delta,\E)
\end{equation}
enhanced-cartesian over $\Unf(\Delta)$. Explicitly, enhanced objects
in the target of \eqref{xi.do} are pairs $\langle [n],E \rangle$,
$[n] \in \Delta$, $E:\Unf((\Delta[n])^o) \to \E$ an enhanced
functor. Say that $E$ is {\em special} if the underlying unenhanced
functor $(\Delta[n])^o \to \E_\ppt$ is special in the sense of
Subsection~\ref{adm.subs}; then by \cite[Lemma 7.4.3.2]{ka.big},
\cite[Lemma 7.4.1.10]{ka.big} and Example~\ref{xi.n}, \eqref{xi.do}
identifies $\Delta^{h\Tot}\E$ with the full subcategory $
\fFun^h(\tw(\Delta)^o|\Delta,\E)_{sp} \subset
\fFun^h(\tw(\Delta)^o|\Delta,\E) $ spanned by pairs $\langle [n],E
\rangle$ with special $E$.

Now take $\E = \D_{[0,1]}(\A)^h$, and note that by
Definition~\ref{adm.def}~\thetag{ii}, the quotient enhanced functor
$h:\Unf(C_{[0,1]}(\A)) \to \D_{[0,1]}(\A)$ induces an enhanced
functor
\begin{equation}\label{q.eq}
Q:\Unf(\wt{C}_{[0,1]}(\A)) \to
\fFun^h(\tw(\Delta)^o|\Delta,\D_{[0,1]}(\A)^h)_{sp}
\cong \Delta^{h\Tot}\D_{[0,1]}(\A)^h,
\end{equation}
enhanced-cartesian over $\Unf(\Delta)$. Both the source and the
target of \eqref{q.eq} are enhanced Segal categories. Moreover,
since any object in $\D_{[0,1]}(\A)$ can be represented by an
admissible functor $\Delta^o \to C_{[0,1]}(\A)$, the underlying
unenhanced functor $Q_\ppt$ is essentially surjective over $[0] \in
\Delta$. Thus to prove that $Q$ induces a desired $2$-equivalence,
it suffices to check that it is enhanced $2$-fully faithful --- that
is, the corresponding square \eqref{ff.1.sq} is
semicartesian. Explicitly, this square is a disjoint union of
squares
\begin{equation}\label{ff.D}
\begin{CD}
\Unf(C^{(2)}(M_\idot,N_\idot)) @>>>
(\Delta^{h\Tot}\D_{[0,1]}(\A)^h)_{[1]}\\
@VVV @VV{s^* \times t^*}V\\
\ppt^h @>{\eps^h(M_\idot) \times \eps^h(N_\idot)}>> \D_{[0,1]}(\A)^h
\times \D_{[0,1]}(\A)^h
\end{CD}
\end{equation}
over all $M_\idot,N_\idot \in C_{[0,1]}(\A)$ representing enhanced
objects in $\D_{[0,1]}(\A)^h$. We need to prove that all the squares
\eqref{ff.D} are semicartesian over all $J \in \Pos$, but since we
can replace $\A$ with $J^o\A$, it suffices to consider the situation
over $\ppt$ --- that is, the underlying unenhanced squares.

Now, the enhanced fiber $(\Delta^{h\Tot}\D_{[0,1]}(\A)^h)_{[1]}$ in
\eqref{ff.D} is the enhanced functor category
$\fFun^h([1],\D_{[0,1]}(\A)^h) \cong \D_{[0,1]}(\Fun([1],\A))^h$, so
its truncation is simply $\D_{[0,1]}(\Fun([1],\A))$. If we denote by
$\D_{[0,1]}(\Fun([1],\A))_{M_\idot,N_\idot}$ the unenhanced fiber of $s^*
\times t^*$ over $M_\idot \times N_\idot \in \D_{[0,1]}(\A) \times
\D_{[0,1]}(\A)$, then \eqref{ff.D} reduces to a functor
\begin{equation}\label{ff.MN}
  C^{(2)}(M_\idot,N_\idot) \to
  \D_{[0,1]}(\Fun([1],\A))_{M_\idot,N_\idot},
\end{equation}
and we need to show that this is an epivalence for any $M_\idot$ and
$N_\idot$.

By definition, the target of the functor \eqref{ff.MN} is the
category of objects $E_\idot \in \D_{[0,1]}(\Fun([1],\A))$ equipped
with isomorphisms $s^*E_\idot \cong M_\idot$, $t^*E_\idot \cong
N_\idot$ in $\D_{[0,1]}(\A)$. Among other things, such an object
$E_\idot$ has homology objects $H_l(E_\idot) \in \Fun([1],\A)$, $l =
0,1$, and these define maps $f_l:H_l(M_\idot) \to H_l(N_\idot)$, so
that we have a projection
$$
\D_{[0,1]}(\Fun([1],\A))_{M_\idot,N_\idot} \to
  \Hom(H_0(M_\idot),H_0(N_\idot)) \times
  \Hom(H_1(M_\idot),H_1(N_\idot))
$$
whose target is a discrete category; the functor \eqref{ff.MN}
intertwines this projection with the projection \eqref{C.MN}. The
source of the functor \eqref{ff.MN} has been described in
Proposition~\ref{cat.prop}. To see \eqref{ff.MN} in terms of
\eqref{spl.3}, assume given a splitting $\wM$ of $(f_1 \oplus
\id) \circ (M_\idot \oplus N_\idot) \circ (\id \oplus (-f_0))$, and
let $\wM_\idot$ be the two-term complex
$$
\begin{CD}
  (M_1 \oplus N_1) @>>> \wM.
\end{CD}
$$
We then have natural maps
\begin{equation}\label{MN.w}
\begin{CD}
M_\idot @<{g}<< \wM_\idot @>{f}>> N_\idot
\end{CD}
\end{equation}
induced by the projections $M_\idot \oplus N_\idot \to
M_\idot,N_\idot$, and $g$ is a quasiisomorphism, while $f$ induces
$f_l$ on $H_l$, $l=0,1$. Treating the arrow $f$ as a two-term
complex in $\Fun([1],\A)$ gives an object $E_\idot \in
C_{[0,1]}(\Fun([1],\A))$ equipped with an isomorphism $t^*E_\idot
\cong N_\idot$, while $g$ provides a quasiisomorphism $s^*E_\idot
\cong M_\idot$. This is where \eqref{ff.MN} sends our splitting. To
finish the proof, we need to show that any object in
$\D_{[0,1]}(\Fun([1],\A)_{M_\idot,N_\idot}$ can be represented in
this way for some splitting, and that any map between such objects
comes from a map of splittings. This is straightforward diagram
chasing that we leave to the reader.
\endproof

\begin{remark}
At the end of the day, enhancements are used two times in our
construction of the $2$-category
$\Delta^h\D_{[0,1]}(\A)^h$. Firstly, we need enhanced functor
categories to obtain the enhanced Segal category
$\Delta^{h\Tot}\D_{[0,1]}(\A)^h$. Secondly, we use semicartesian
products to cut it down to $\Delta^h\D_{[0,1]}(\A)^h$ via
\eqref{red.h.sq}. The first instance is not too serious ---
effectively, at least if we only want the truncation
$\Delta^{h\Tot}\D_{[0,1]}(\A)^h_{\ppt}$, we might as well take the
simplicial replacement $\Delta^{\Tot}C_{[0,1]}(\A)$ and localize it
with respect to quasiisomorphisms vertical over $\Delta$. This
produces a fibration with fibers $\D_{[0,1]}(\Fun([n],\A))$, just as
expected. However, the ability to take semicartesian products is
absolutely crucial. Sometimes, for specially chosen $M_\idot$ and
$N_\idot$, the functor \eqref{ff.MN} is an equivalence, and were it
always the case, it would have been enough to take the cartesian
product. However, for general $M_\idot$ and $N_\idot$, \eqref{ff.MN}
is really only an epivalence --- for example, if $M_\idot =
N_\idot$, and we consider the locus where $f_0=\id$, $f_1=\id$, then
the target of \eqref{ff.MN} is discrete, while the source is a
groupoid with non-trivial $\pi_1$ (isomorphic to
$\Hom(H_0(M_\idot),H_1(N_\idot))$). Thus we do really need to use
the general machinery, and the black box result of \cite[Lemma
  7.3.3.7]{ka.big}.
\end{remark}

\begin{remark}
One might wonder what happens when one looks at the full derived
category $\D(\A)^h$ rather than just $\D_{[0,1]}(\A)^h$. In this
case, we still have an enhanced $2$-category --- in fact, the
enhanced Yoneda pairing of \cite[Remark 7.4.6.3]{ka.big} does
provide the homotopy type of maps $\Hhom(M,N)$ for any objects $M,N
\in \D(\A)$, and \eqref{red.h.sq} still works and provides an
enhanced $2$-category refining $\D(\A)$. However, it is no longer
$1$-truncated. Its $1$-truncation is a $2$-category in the usual
unenhanced sense, but we do not know whether it admits a reasonably
concise explicit construction in the spirit of
Subsection~\ref{2.co.subs}.
\end{remark}

{\small\noindent
Affiliations:
\begin{enumerate}
\renewcommand{\labelenumi}{\arabic{enumi}.}
\item Steklov Mathematics Institute (main affiliation).
\item National Research University Higher School of Economics.
\end{enumerate}}

{\small\noindent
{\em E-mail address\/}: {\tt kaledin@mi-ras.ru}
}

\end{document}